\documentclass[10pt]{article}

\usepackage{amscd,amsmath, amssymb, fancyhdr}

\usepackage[backref=page]{hyperref}
\renewcommand*{\backrefalt}[4]{%
	\ifcase #1 (Not cited.)%
	\or        (Cited on page~#2.)%
	\else      (Cited on pages~#2.)%
	\fi}

\hypersetup{
	colorlinks   = true,
	citecolor    = magenta,
	linkcolor    = blue,
	urlcolor     = magenta	
}

\numberwithin{equation}{section}

\newcommand{\version}{version 2.1,\ \ Dec. 27, 2025}

\def\eqref#1{(\ref{#1})}

\newcommand{\Z}{{\Bbb Z}}
\def\C{{\Bbb C}}

\newcommand{\Q}{{\Bbb Q}}

\def\1{\sqrt{-1}\:}

\newcommand{\cntrct}                
{\hspace{2pt}\raisebox{1pt}{\text{$\lrcorner$}}\hspace{2pt}}

\renewcommand{\phi}{\varphi}
\renewcommand{\epsilon}{\varepsilon}

\renewcommand{\min}{{\it min}}

\newcommand{\Aut}{\operatorname{Aut}}
\newcommand{\Alb}{\operatorname{Alb}}

\renewcommand{\gcd}{\operatorname{ \text {\it \sffamily gcd}}}

\newcounter{Mycounter}[section]
\newcounter{lemma}[section]
\newcounter{claim}[section]
\renewcommand{\theclaim}{{Claim \thesection.\arabic{claim}}}
\newcommand{\claim}{%
    \setcounter{claim}{\value{Mycounter}}
    \refstepcounter{claim}
    \stepcounter{Mycounter}
    {\noindent \bf \theclaim:\ }}

\newcounter{sublemma}[section]
\newcounter{corollary}[section]
\renewcommand{\thecorollary}{{Corollary \thesection.\arabic{corollary}}}
\newcommand{\corollary}{%
    \setcounter{corollary}{\value{Mycounter}}
    \refstepcounter{corollary}
    \stepcounter{Mycounter}
    {\noindent \bf \thecorollary:\ }}

\newcounter{theorem}[section]
\renewcommand{\thetheorem}{{Theorem \thesection.\arabic{theorem}}}
\newcommand{\theorem}{%
    \setcounter{theorem}{\value{Mycounter}}
    \refstepcounter{theorem}
    \stepcounter{Mycounter}
    {\noindent \bf \thetheorem:\ }}

\newcounter{conjecture}[section]
\newcounter{proposition}[section]
\newcounter{definition}[section]
\renewcommand{\thedefinition}
      {{Definition~\thesection.\arabic{definition}}}
\newcommand{\definition}{%
    \setcounter{definition}{\value{Mycounter}}
    \refstepcounter{definition}
    \stepcounter{Mycounter}
    {\noindent \bf \thedefinition:\ }}

\newcounter{example}[section]
\newcounter{remark}[section]
\renewcommand{\theremark}{{Remark \thesection.\arabic{remark}}}
\newcommand{\remark}{%
    \setcounter{remark}{\value{Mycounter}}
    \refstepcounter{remark}
    \stepcounter{Mycounter}
    {\noindent \bf \theremark:\ }}

\newcounter{problem}[section]
\newcounter{question}[section]
\newcommand{\pstep}{{\bf Proof. Step 1:\ }}
\newcommand{\proof}{{\bf Proof:\ }}

\makeatletter

\@addtoreset{equation}{section} \@addtoreset{footnote}{section}
\makeatother

\def\blacksquare{\hbox{\vrule width 5pt height 5pt depth 0pt}}
\def\endproof{\blacksquare}

\begin{document}

\begin{center}
{\LARGE\bf
Minimal multiplicity   of fiber\\[2mm] components  in abelian fibrations\\[4mm]
}

Frederic Campana, 
Ljudmila Kamenova\footnote{Partially supported 
by award SFI-MPS-TSM-00013537 
from the Simons Foundation International.}, 
Misha Verbitsky\footnote{Partially supported 
by FAPERJ SEI-260003/000410/2023 and CNPq - Process 310952/2
021-2. 

}

\end{center}

{\small \hspace{0.10\linewidth}
\begin{minipage}[t]{0.85\linewidth}
{\bf Abstract.} An abelian fibration is a 
proper projective surjective map of complex varieties
with general fiber an abelian variety.
Consider a multiple fiber of an
abelian fibration, and let $\mu_1, ..., \mu_k$
be the multiplicities of its irreducible
components. We prove that the minimum of $\mu_i$
is equal to their greatest common divisor
$\gcd(\mu_1, ..., \mu_k)$. 
\end{minipage}
}

\tableofcontents


\section{Introduction} 



Let $f:X\to B$ be a proper fibration from a complex manifold $X$ onto a 
smooth complex curve $B$, and $b\in B$ be a point. 
Let $X_b:=f^*(b):=\sum_i m_iF_i$ be the scheme theoretic (divisorial, here) 
fiber of $f$, the positive integers $m_i$ being the ``multiplicities'' of 
the components $F_i$ of $X_b$. The ``multiplicity'' of $X_b$ is classically 
defined as being the greatest common divisor of the numbers $m_i$. This is, in 
particular, the definition used by Kodaira and Shafarevich in their 
classification of degenerate fibers in elliptic fibrations of surfaces
(\cite{_BHPV:04_}). 
For higher dimensional bases $B$, one looks at the generic fibers over the 
components of the discriminant divisor of $f$ in $B$. 

In \cite{_C'_}, however, a different definition is used: one does not 
consider there the $\gcd$ of the $m_i$, but their infimum. This second 
definition is dictated by the need to characterize what is called there the 
``Bogomolov sheaves'' on projective or compact K\"ahler varieties $X$, 
which are the rank-one subsheaves $L\subset \Omega^p_X, p>0$ such that 
$\kappa(X,L)=p$, the top Kodaira dimension, after Bogomolov's bound
(\cite{_Bogomolov:1979_,_CCR:Specialness_}). 
These sheaves, analogs of ``big'' canonical bundles, indeed play a central 
r\^ole in bimeromorphic classification (see \cite{_C'_} for more details). 

It is thus important to know when the two definitions agree (for locally 
projective or K\"ahler proper fibrations). This is the case when the smooth 
fibers of $f$ are either rationally connected (by
\cite{_GHS_} or \cite[main theorem]{_dJ_S:rational_}), 
 or are compact complex tori (in particular: Abelian varieties), which 
is the topic of the present text. In particular, the two definitions agree 
for the elliptic surfaces of Kodaira/Shafarevich's classification.

On the other hand, examples of fibrations with fibers higher genus curves 
for which the two versions differ are found in \cite{_W_}, \cite{_S_}, 
\cite{_C"_}. A remarkable example is given in \cite{_L_}, where the smooth 
fibers are Enriques surfaces, but some singular fibre has components of multiplicities 2 and 3, the base being the projective line, 
with consequences in hyperbolic and arithmetic 
geometries derived in \cite{_BCJW_}. 

\hfill

In the Appendix \ref{_Torsor:Section_} 
we introduce the abelian torsor structure in the 
settings of complex torus fibrations, as in 
\cite{_BKV:Sections_, _KV:torsors_}. We use the 
averaging formalism applied to multisections of an abelian fibration in order 
to obtain the following theorem over a $1$-dimensional base.

\hfill

\theorem (\ref{_gcd_1_then_section_Theorem_}) 
Let $\pi:\; M \to B$ be an abelian fibration, that is,
a proper projective surjective map of complex varieties with general
fiber an abelian variety. Assume that all points of $B$
are regular values of $\pi$ except $x_0\in B$,
and $\dim_\C B=1$. Let $F_1, ..., F_k$ be the 
irreducible components of $\pi^{-1}(x_0)$,
and $\mu_1, ..., \mu_k$ their multiplicities.
Assume that the greatest common divisor of $\mu_i$
is 1. Then $\pi$ admits a local section.

\hfill

When $\dim_\C M=2$,
\ref{_gcd_1_then_section_Theorem_} is a consequence
of Kodaira's classification of elliptic fibrations,
\cite{_BHPV:04_}. After we posted this paper on the arXiv, we became aware
that \ref{_gcd_1_then_section_Theorem_} was already
proved in an unpublished paper of S. Lu, 
\cite{_Lu:preprint_1999_}.

\hfill

In Section \ref{Section3}, we generalize these results to 
higher-dimensional base. 

\hfill

\begin{theorem} (\ref{prop}) Let $f:M\to B$ be a proper K\"ahler
fibration with $M$ smooth, $B$ normal, such that its generic fibers are 
smooth complex tori. Assume that $f$ has equidimensional fibers.
Let $b\in B$ be any point and let $F_1, ..., F_k$ be the 
irreducible components of $f^{-1}(b)$,
and $\mu_1, ..., \mu_k$ their multiplicities.
Assume that the greatest common divisor of $\mu_i$ is 1. 
Then for any $b\in B$, $f$ has a meromorphic section above some open 
neighborhood of $b$ in $B$. 
\end{theorem}

\hfill

\begin{corollary} (\ref{cor}) In the situation of \ref{prop}, 
assume that for some $b\in B$, the fiber is (in the sense of \cite{_C_}) 
$M_b=\sum_im_iF_i$, with $F_i$ being the irreducible components. 
If $gcd_i(\{m_i\})=1$, then $m_i=1$ for some $i$. 
\end{corollary}

\hfill 

\remark
In assumptions of \ref{prop},
it is essential that the fibration is (at least locally) K\"ahler:
one of the key arguments in Section \ref{Section3} uses the compactness
of the Barlet spaces.
However, in the settings of hyperk\"ahler geometry, the first-named author 
proved that Lagrangian fibrations on hyperk\"ahler manifolds are 
not only K\"ahler, but also locally projective (\cite{Campana}).  


\section{Averaging multisections in a torus fibration over a disk}


In Section \ref{_Torsor:Section_}, we define the torsors and
described the averaging formalism (\cite{_BKV:Sections_}), 
which is a general
method of constructing a section of a torsor bundle,
starting with a (signed) multisection of total multiplicity 1.
Applying this construction to an abelian fibration, we obtain
the following theorem.

\hfill

\theorem\label{_gcd_1_then_section_Theorem_}
Let $\pi:\; M \to B$ be a torus fibration, that is,
a proper map of complex varieties with general
fiber a compact torus. Assume that all points of $B$
are regular values of $\pi$ except $x_0\in B$,
and $\dim_\C B=1$. Let $F_1, ..., F_k$ be the 
irreducible components of $\pi^{-1}(x_0)$,
and $\mu_1, ..., \mu_k$ their multiplicities.
Assume that the greatest common divisor of $\mu_i$
is 1. Then $\pi$ admits a local section, and
$\min_i \mu_i=1$.

\hfill

\pstep
Around a smooth point $m\in F_i$,
the projection maps can be written in local coordinates
as $(z_1, z_2, ..., z_n) \mapsto z_1^{\mu_i}$. 
Let ${\cal E}_{\mu_i}$ be the set of all roots of unity of degree $\mu_i$.
Then the image of the 
multivalued map $t \to ({\cal E}_{\mu_i}t^{\mu_i^{-1}}, 0,..., 0)$
is a smooth curve in $M$, and its projection 
to a neighbourhood $U$ of $x\in B$ is a $\mu_i$-sheeted ramified cover.

\hfill

{\bf Step 2:}
Consider such a multisection $Z_i$ for each $F_i$, and
let $a_1,..., a_k$ be a set of integers such that
$\sum_i a_i \mu_i =1$. We consider $\sum a_i Z_i$
as a multisection with multiplicites determined by 
$a_i$, with total multiplicity 1. 
Applying the averaging construction 
(\ref{_averaging_Torsor_Claim_}), we obtain a meromorphic
map $\phi:\; U\backslash x_0 \to M$, which is holomorphic
and Lipschitz 
by construction. By Riemann removable singularity theorem, 
this map extends to the origin $x_0\in U$.
\endproof

\hfill

\corollary
Let $\pi:\; M \to B$ be an torus fibration, that is, 
a proper surjective map of complex varieties with general 
fiber a compact torus. Assume that all points of $B$
are regular values of $\pi$ except $x_0\in B$,
and $\dim_\C B=1$. Let $F_1, ..., F_k$ be
the irreducible components of $\pi^{-1}(x_0)$,
and $\mu_1, ..., \mu_k$ their multiplicities.
Then $\min \mu_i = \gcd(\mu_i)$,
where $\gcd$ denotes the greatest common divisor.

\hfill

\proof
Passing to a ramified covering of $X$ of order $\gcd(\mu_i)$,
we arrive to an torus fibration with the same fibers with multiplicities
$\mu_i'= \frac{\mu_i}{\gcd(\mu_i)}$, hence $\gcd(\mu_i')=1$.
By \ref{_gcd_1_then_section_Theorem_}, it has a section.
This section gives a multisection of order $\gcd(\mu_i)$ of $\mu$.
\endproof


\section{Higher dimensional bases} \label{Section3}


The aim of this section is to prove the following results generalizing 
\ref{_gcd_1_then_section_Theorem_} from a $1$-dimensional base to a base of 
any dimension. First we give a proof of the corollary, and then of the 
theorem below. 

\hfill

\begin{theorem}\label{prop} Let $f:M\to B$ be a proper K\"ahler
fibration between normal complex spaces, $M$ smooth, 
such that its generic fibers are 
smooth compact complex tori. Assume that $f$ has 
equidimensional fibers. Let $b\in B$ be any point 
and let $F_1, ..., F_k$ be the 
irreducible components of $f^{-1}(b)$,
and $\mu_1, ..., \mu_k$ their multiplicities.
Assume that the greatest common divisor of $\mu_i$
is 1.
Then $f$ has a meromorphic section above some open 
neighborhood of $b$ in $B$.  \end{theorem}

\hfill

\begin{corollary} \label{cor} In the situation of \ref{prop}, 
assume that for some $b\in B$, the fiber is (in the sense of \cite{_C_}) 
$M_b=\sum_im_iF_i$, with $F_i$ being the irreducible components. 
If $gcd_i(\{m_i\})=1$, then $\min_i (\{m_i\})=1$. 
\end{corollary}

\hfill 

{\bf Proof of \ref{cor}, assuming \ref{prop}:} Let $s:B\dasharrow M$ be a 
meromorphic section of $f$ defined near $b$. Let $B'$ be the closure
of the image of the meromorphic section $B \dasharrow M$. Let $m:B'\to B$ be the proper 
bimeromorphic map such that $s':=s\circ m:B'\to M':=M\times_B B'$ is a 
holomorphic section of $f':M'\to B'$ obtained from $f$ by the base-change 
$m:B'\to B$. 
For any $b'\in B'$ such that $m(b')=b$, the component of 
$M'_{b'}$ containing $s'(b')$ has multiplicity $1$ in 
$M'_{b'}=M_b\times_{\{b\}}\{b'\}$. Thus, so has the corresponding component 
of $M_b$: since the existence of a local holomorphic section implies that 
the multiplicity of the corresponding component of the fiber has multiplicity one. 
\endproof

\hfill

{\bf Proof of \ref{prop}:} 
Let $B^0\subset B$ be the Zariski open subset 
over which 
the fibers are smooth compact tori. Fix $b\in B^0$, let $A_b=\Aut^0(M_b)$ 
be the translation group of the fiber $M_b$. Clearly, $A_b$ is a torus,
isomorphic to $M_b$.
If $u_i\in \Bbb Z,i=1,\dots,r$ is such that $\sum_iu_i=1$ and if 
$p_i\in M_b,i=1,\dots, r$, the sum $\sum_iu_i\cdot p_i=s\in M_b$ is well-defined, 
independently of the choice of an origin in $M_b$ 
(\ref{_averaging_Torsor_Claim_}).  

Recall (\cite{_C_}, chap IV, Th\'eor\`eme 3.4.1, Corollaire 3.4.2) that $f$ being equidimensional, with $B$ normal, 
the components of the fibers of $f$ can be given in a unique way 
multiplicities (equal to $1$ for those above $B^0$) in such a way that the 
fibers of $f$ form an analytic family of cycles in the Barlet sense, 
such that the map $\beta: B\to C(M/B)$ sending $M_b$ to its representative 
point in $C(M/B)$, the Chow/Barlet space of relative $f$-cycles of $M$, 
is analytic. Moreover, if $C:=\beta(B)\subset C(M/B)$, and if 
$\Gamma_C\subset C\times M$ is the incidence graph of the universal family 
of cycles parametrised by $C$, then 
$(f\times id_X)(M)=(\beta\times id_X)^{-1}(\Gamma_C)$: by its very definition, 
the graph of $f$ is obtained from $\Gamma_C$ by the base-change $\beta$.

Let now $b'\in B\setminus B^0$, and $M_{b'}=\sum_im_i\cdot
F_i$ be its $f$-fiber. 
By smoothness of $M$, through the generic point of each $F_i$ we have a 
local multisection $s_i,i=1,\dots, r$ of $f$ which has degree $m_i$ over $B$. 
Denote with $p_{i,j}(b)$ its points over $b$.  By hypothesis, the gcd of the 
numbers $\{m_i\}$ 
is one, and Bezout lemma
implies the existence of integers $v_i$ 
such that $\sum_i u_i:=\sum_iv_i\cdot m_i=1, u_i:=v_i\cdot
m_i$.   By averaging (\ref{_averaging_Torsor_Claim_}), 
for each $b\in B^0$ near $b'$, a well-defined unique point  
$s(b):=\sum_{i,j} u_i\cdot p_{i,j}(b)\in M_b$ exists, which depends holomorphically on $b$. 

We shall show, using the K\"ahlerianity assumption, that this section of $f$ 
over $B^0$ near $b'$ extends meromorphically to $B$ near $b'$. 
More precisely, we will show that the addition (and substraction) in the generic smooth fibers of $f$ extend 
meromorphically to $M\times_B M$. Any section obtained by averaging over the locus $B^0$ of smooth fibers 
of $f$ in $B$ thus extends meromorphically to $B$.

This extension of addition/substraction is based on the properness of the irreducible components of 
the relative Chow/Barlet space $C(M/B)$ (see \cite{_C_} and the references 
there). 
Let $B^0\subset B$ be the Zariski open subset over which the fibers are 
smooth compact tori, and let $M^0:=f^{-1}(B^0)$. 

Over $B^0$ one can define the unit component $\Aut^0(M^0/B^0)\to B^0$ of the 
group of translations in the fibers of $f$ as the component of 
$C(M^0\times_{B^0}M^0)\to B^0$ consisting of the graphs of these translations. 
Over $B^0$, the incidence graph 
$\Phi^0\subset (\Aut^0(M^0/B^0)\times_{B^0}(M^0\times_{B^0}M^0)$ 
of the family of relative translations in the fibers of $f$ is also nothing 
but the graph of the addition map $+:\Aut^0(M^0/B^0)\times M^0\to M^0$ over 
$B^0$ sending $(t_b,x_b)\to t_b+x_b$, the substraction map being defined as 
exchanging the $M^0/B^0$ factors. Since over $B^0$ these two maps just consist 
in exchanging the suitable factors, which obviously extends meromorphically over $B$, 
they extend meromorphically over $B$ to maps: $\Aut^0(M/B)\times_BM_b$, where $Aut^0(M/B)$ is the Zariski 
closure in $C(M/B)$ of the unit component $\Aut^0(M^0/B^0)$. 
If $u_i\in \Bbb Z,i=1,\dots r$ are such that $\sum_iu_i=1$, the 
well-defined map over $B^0$: $M^0\times_{B^0}\times \dots _{B^0}M^0\to M_0$ 
(r-factors) sending $(x_{1,b},\dots x_{r,b})\to
\sum_iu_i\cdot x_{i,b}$ 
thus extends meromorphically to $M_B\times\dots_BM$ over $B$, 
which shows the claim.
\endproof

\hfill 

\begin{remark} \label{rem}
1. Constructions similar to the one given by Ulf Persson in \cite[p.83]{_P_}, 
show that the relative K\"ahlerianity of $f$ is essential for the existence 
of this meromorphic section.

2. The hypothesis on the $\gcd$ of the $m_i$ being 
equal to $1$ does not follow from the same property in codimension $1$, 
even if $M$ and $B$ are smooth, $f$ projective, and if the fibers outside 
$b\in B$ are reduced, as shown by the example of I. Hellmann 
(\cite[Section 2.2]{_Hellmann:cone_}). 
\end{remark}

\section{Appendix: Torsor fibrations} \label{_Torsor:Section_}


\subsection{Torsor structure in complex torus fibrations}

By definition, a torsor over a group
$G$ is a set where $G$ acts freely and transitively.
However, the torsor structure can be defined
without referring to an {\em a priori} given group structure as follows.

\hfill

Following \cite{_KV:torsors_}, we define an 
abelian torsor structure as a kind of algebraic structure,
without first defining the group acting on this torsor.
This simplifies some of the definitions when the
torsor structure is defined in families, and we
do not have the corresponding family of abelian
groups on hand.

\hfill

In the sequel we will make a distinction between the compact
torus $T$ and the same torus equipped with a group structure,
which we denote $\Alb(T)$. This group structure is clearly
unique up to a parallel transport, that is, determined
by the choice of 0. The torus $\Alb(T)$ is identified
with the abelian variety of $T$.

\hfill

\definition
{\bf An abelian torsor} is a set $S$ equipped
with a map $\mu:\; S\times S \times S \to S$ such that
\begin{description}
\item[(i)] for
any $e\in S$, the map $x \cdot y := \mu(x, y, e)$
defines a commutative and associative product on $S$
with $e$ an identity element 
\item[(ii)] 
 for any fixed $y \in S$,
the map $x\mapsto x\cdot y$ is bijective.
\end{description}

\remark
From the axioms of an abelian
torsor, it is clear that the operation ``$\cdot$'' defines a structure 
of an abelian group on $S$.

\hfill

\remark
The identity element axiom means that 
$\mu(x,y,y)=x$ for any $x$ and $y$. Given a group
structure on $S$, the torsor operation
$\mu$ can be defined  as $\mu(x, y,z)= x+y-z$.

\hfill

\claim\label{_fibrations_torsor_Claim_}
Let $\pi:\; M \to B$ be 
a surjective holomorphic fibration whose general fiber is a 
compact complex torus.
Then the fibers $F_x:= \pi^{-1}(x)$ of $\pi$ are
equipped with a holomorphic torsor structure
which depends on $x\in B$ holomorphically.

\hfill

\proof
It is not hard to see that any 
holomorphic map of a compact complex torus
is affine, that is, obtained as a composition
of a parallel transport and a group homomorphism.
Therefore, the holomorphic torsor structure on a
compact torus $T$ is unique; indeed, the 
connected component of the group of automorphisms
of $T$ is $\Alb(T)$. Taking a local section,
we can identify each fiber $F_x$ with its
Albanese manifold $\Alb(F_x)$ holomorphically,
and this turns $F_x$ into an  $\Alb(F_x)$-torsor. 
The dependence on $x$ is holomorphic by construction,
and it is independent from the choice of a local section
by uniqueness.
\endproof

\subsection{The averaging formalism}

An affine space is, by definition, a torsor over a vector space.
Consider an affine space $A$ over a vector space $V$. 
After choosing the origin $a\in A$, we can identify
$A$ and $V$. For any collection of points $x_1, ..., x_n\in A$,
their mass center $a + \frac 1 n \sum_i(x_i-a)$ 
is well defined and independent from the choice of the origin.

\hfill

This construction works for torsors over infinitely
divisible, torsion-free abelian groups (that is, vector spaces
over $\Q$). This assumption is essential, because if the group has 
torsion, the mass center becomes ambiguous.

\hfill

However, a version of this construction
survives. It was used in \cite{_BKV:Sections_}
to produce a rational section of a Lagrangian fibration
in a hyperk\"ahler manifold. If we take a weighted
sum of $m$ points with positive sign and $(m-1)$ points
with negative sign, the sum is independent from the
choice of an origin.
We give an abstract version of this construction here, see 
\cite[Claim 3.6]{_BKV:Sections_} for its proof.

\hfill

\claim\label{_averaging_Torsor_Claim_}
Let $X$ be a torsor over $G$.
Choosing an origin $e\in X$, we can identify
$X$ with $G$. Then, for any $m\in \Z^{>0}$ there is a natural
map $X^m\times X^{m-1} \to X$
taking $(x_1,..., x_m, y_1, ..., y_{m-1})$
to $e+ \sum_{i=1}^{m} (x_i-e) - \sum_{i=1}^{m-1} (y_i-e)$,
which is independent from the choice of $e$.
\endproof

\hfill

\hfill

{\bf Acknowledgments.} 
L. K. and M. V. are grateful to J. Koll\'ar 
for his insightful comments and email exchanges, 
and to F. Bogomolov, whose ideas were of great 
influence in this paper. We express our gratitude to
J. Koll\`ar for noticing that, in the first 
version of this paper, the gcd hypothesis 
had been forgotten in the statement of Theorem \ref{prop}. 
We are also grateful to
S. Lu for communicating to us the results of
\cite{_Lu:preprint_1999_}.

The present text was written during a visit of F. C. at IMPA in Rio, and he 
thanks this institute and the local branch of CNRS for their support, making 
this visit possible. 

\hfill

\small
\noindent {\sc Frederic Campana\\
D\'epartement de math\'ematiques.
Universit\'e Nancy 1.\\
64, Boulevard des Aiguillettes.\\
54506 Vandoeuvre-l\`es-Nancy C\' edex. France.}\\
\tt frederic.campana@univ-lorraine.fr
\\

\noindent {\sc Ljudmila Kamenova\\
Department of Mathematics, 3-115 \\
Stony Brook University \\
Stony Brook, NY 11794-3651, USA,} \\
\tt kamenova@math.sunysb.edu
\\

\noindent {\sc Misha Verbitsky\\
            {\sc Instituto Nacional de Matem\'atica Pura e
              Aplicada (IMPA) \\ Estrada Dona Castorina, 110\\
Jardim Bot\^anico, CEP 22460-320\\
Rio de Janeiro, RJ - Brasil}\\
\tt verbit@impa.br
}

\end{document}